\documentclass[12pt,twoside]{amsart}

\usepackage[dvipsnames]{xcolor}
\usepackage[normalem]{ulem}
\usepackage{hyperref}
\hypersetup{colorlinks=true,
allcolors=Blue}
\usepackage{amsthm, amsmath, amscd, amssymb,centernot,txfonts}
\usepackage{tikz-cd}
\usepackage{lipsum}
\usepackage{amsfonts}
\usepackage{listings}
\usepackage{multicol}
\usepackage{array}
\usepackage{multirow}
\usepackage{mathtools}
\usepackage{tabularx,ragged2e,booktabs,caption}

\usepackage{mathrsfs}
\usepackage[all]{xy}
\usepackage[left=2cm,top=2.5cm,bottom=3cm,right=3cm]{geometry}
\usepackage{dirtytalk}

\newcolumntype{M}[1]{>{\centering\arraybackslash}m{#1}}
\theoremstyle{plain}
\numberwithin{equation}{section}
\newtheorem{theorem}{Theorem}[section]
\newtheorem{proposition}[theorem]{Proposition}
\newtheorem{lemma}[theorem]{Lemma}
\newtheorem{corollary}[theorem]{Corollary}

\newtheorem{definition}[theorem]{Definition}
\newtheorem{remark}[theorem]{Remark}

\newtheorem{observation}[theorem]{Observation}

\newcommand*{\QEDB}{\hfill\ensuremath{\square}}
\begin{document}

\title[Remarks on projective normality for calabi-yau and hyperk\"ahler varieties]{Remarks on projective normality for certain Calabi-Yau and hyperk\"ahler varieties}

\author[Mukherjee]{Jayan Mukherjee}
\address{Department of Mathematics, University of Kansas, Lawrence, KS 66045}
\email{j899m889@ku.edu}

\author[Raychaudhury]{Debaditya Raychaudhury}
\address{Department of Mathematics, University of Kansas, Lawrence, KS 66045}
\email{debaditya@ku.edu}

\maketitle

\begin{abstract}
We prove some results on effective very ampleness and projective normality for some varieties with trivial canonical bundle. In the first part we prove an effective projective normality result for an ample line bundle on regular smooth four-folds with trivial canonical bundle. More precisely we show that for a regular smooth fourfold with trivial canonical bundle, $A^{\otimes 15}$ is projectively normal for $A$ ample. In the second part we emphasize on the projective normality of multiples of ample and globally generated line bundles on certain classes of known examples (upto deformation) of projective hyperk\"ahler varieties. As a corollary we show that excepting two extremal cases in dimensions $4$ and $6$, a general curve section of any ample and globally generated linear system on the above mentioned examples is non-hyperelliptic.   
\end{abstract}

\section*{Introduction}
In this article, we prove some effective projective normality results for multiples of line bundles with some positivity conditions on Calabi-Yau and hyperk\"ahler varieties. These are varieties with trivial canonical bundles that are higher dimensional analogues of $K3$ surfaces. We recall some results on syzygies for varieties with trivial canonical bundle below and state our main result for regular four-folds with trivial canonical bundles. Next we will concentrate on a special class of even dimensional such varieties (regular with trivial canonical bundle), namely hyperk\"ahler varieties. 
\subsection{Regular varieties with trivial canonical bundle.}
Geometry of linear series on varieties $X$ with trivial canonical bundle $K_X$ (i.e. $K_X=0$) is a topic that has motivated a lot of research. The question of what multiple of an ample line bundle on a $K$-trivial variety is very ample and projectively normal was extensively studied by many mathematicians including Gallego, Oguiso, Peternell, Purnaprajna, Saint-Donat (see \cite{GP1}, \cite{OP}, \cite{SD}). Recall that a $K3$ surface is defined as a smooth projective surface $S$ with $K_S = 0$ and $H^1(\mathscr{O}_S)=0$. The following theorem is due to Saint-Donat (see \cite{SD}). \par
\phantomsection\label{A}\begin{theorem}
Let $S$ be a smooth projective $K3$ surface and let $B$ be an ample line bundle on $S$. Then $B^{\otimes n}$ is very ample (in fact projectively normal) for $n\geq 3$.
\end{theorem}
 Gallego and Purnaprajna generalized Saint-Donat's result on projective normality for smooth, projective, regular ($H^1(\mathscr{O}_X)=0$) three-fold $X$ with trivial canonical bundle (see \cite{GP1}, Corollary 1.10). They showed that for an ample line bundle $B$ on $X$, $B^{\otimes n}$ is projectively normal for $n=8$ and $n\geq 10$. Moreover, if $B^3>1$ then $B^{\otimes n}$ is projectively normal if $n=6$ and $n\geq 8$.
In order to prove this theorem, Gallego and Purnaprajna studies the case when a regular, three-fold with trivial canonical bundle maps onto a variety of minimal degree by a complete linear series of an ample and globally generated line bundle. They gave a classification theorem for such situations and proved the projective normality result as a consequence by applying a theorem of Green (see Theorem \hyperref[green]{1.7}). Varieties that appear as covers of varieties of minimal degree play an important role in the geometry of algebraic varieties. They are extremal cases in a variety of situations from algebraic curves to higher dimensional varieties (see \cite{GP1}, \cite{G1}, \cite{G2}, \cite{Hor}).\\
\indent In this article, we prove an analogue of the classification theorem of Gallego and Purnaprajna where we study the situation when a smooth regular  four-fold $X$ with trivial canonical bundle maps to a variety of minimal degree by the complete linear system of an ample and globally generated line bundle $B$ (see Theorem \hyperref[class]{2.3}) and provide upper bounds to the degree of such morphisms. As a consequence of the classification theorem, Fujita's base point freeness conjecture that has been proved in dimension four by Kawamata (see \cite{Ka}) and Green's theorem (c.f. Theorem \hyperref[green]{1.7}), we prove the following effective projective normality result (see Theorem \hyperref[pthma]{2.4}) on regular fourfolds with trivial canonical bundle.\par 

\noindent {\bf Theorem A.} \phantomsection\label{oa}\textit{Let $X$ be a four-fold with trivial canonical bundle and let $A$ be an ample line bundle on $X$. then, 
\begin{itemize}
    \item[(1)] $nA$ is very ample and embeds $X$ as a projectively normal variety for $n \geq 16$.
    \item[(2)] If $H^{1}(\mathscr{O}_{X})=0$ then $nA$ is very ample and embeds $X$ as a projectively normal variety for $n \geq 15$.
\end{itemize}}\par     

This theorem can be thought of as a higher dimensional analogue of Theorem \hyperref[A]{0.1} and Corollary 1.10, \cite{GP1}. Note that standard methods of Castelnuovo-Mumford regularity (see Lemma \hyperref[cm]{1.6}) and Theorem 1.3, \cite{GP1} yields in the situation above that $A^{\otimes n}$ satisfies projective normality for $n\geq 21$.\\
\subsection{Calabi-Yau and Hyperk\"ahler varieties.} The definition of $K3$ surface, i.e. a smooth regular surface $S$ with trivial canonical bundle is equivalent to having a holomorphic symplectic form on $S$. However in higher dimensions these two notions do not coincide which is clear from the fact that existence of a holomorphic symplectic form on a K\"ahler manifold demands that its dimension is even whereas there are examples of smooth projective algebraic varieties in odd dimensions with trivial canonical bundle and $H^1(\mathscr{O}_X)=0$, for example smooth hypersurfaces of degree $n+1$ in $\mathbb{P}^n$. So essentially we can have two different kinds of generalizations of a $K3$ surface. We give the precise definitions of Calabi-Yau and hyperk\"ahler varieties below. 
\begin{definition}
\phantomsection\label{0.1} A compact K\"ahler manifold $M$ of dimension $n\geq 3$ is called Calabi-Yau if it has trivial canonical bundle and the hodge numbers $h^{p,0}(M)$ vanish for all $0<p<n$.
\end{definition} 
With this definition, Calabi-Yau manifolds are neccessarily projective. However the following definition of hyperk\"ahler manifolds does not imply projectivity in general.
\begin{definition}
\phantomsection\label{0.2}A compact K\"ahler manifold $M$ is called hyperk\"ahler if it is simply connected and its space of global holomorphic two forms is spanned by a symplectic form.
\end{definition}
The decomposition theorem of Bogomolov (see \cite{Bo}) says, any complex manifold with trivial first Chern class admits a finite \'etale cover isomorphic to a product of complex tori, Calabi-Yau manifolds and hyperk\"ahler manifolds. Thus, these spaces can be thought as the ``building blocks'' for manifolds with trivial first Chern class.\\
\indent The theorem of Saint-Donat that we stated (see Theorem \hyperref[A]{0.1}) deals with ample line bundles. In the same paper (see \cite{SD}), Saint-Donat proves the following theorem for ample and globally generated line bundles on $K3$ surfaces. \par
\begin{theorem}\phantomsection\label{SDR}
Let $S$ be a smooth projective $K3$ surface and let $B$ be an ample and base point free line bundle on $S$. Then, \begin{itemize}
    \item[(1)] $B^{\otimes 2}$ is very ample and $|B^{\otimes 2}|$ embeds $S$ as a projectively normal variety unless the morphism given by the complete linear system $|B|$ maps $S$, $2:1$ onto $\mathbb{P}^2$.
    \item[(2)] $B$ is very ample and $|B|$ embeds $S$ as a projectively normal variety unless the morphism given by the complete linear system $|B|$ maps $S$, $2:1$ onto $\mathbb{P}^2$ or to a variety of minimal degree.
\end{itemize}
\end{theorem} 
Gallego and Purnaprajna provides a generalization of this theorem for regular three-folds $X$ with trivial canonical bundle in \cite{GP1} in which they proved for an ample and globally generated line bundle $B$, $B^{\otimes 3}$ is projectively normal unless the morphism $\varphi_B$ induced by the complete linear series of $B$ maps $X$, 2-1 onto $\mathbb{P}^3$. Moreover, they showed that $B^{\otimes 2}$ is projectively normal unless $\varphi_{B}$ maps $X$, 2-1 onto $\mathbb{P}^3$ or to a variety of minimal degree $\geq 2$. They also proved that $B^{\otimes 4}$ is projectively normal on smooth, projective, regular, four-folds with trivial canonical bundle when the morphism $\varphi_B$ associated to the complete linear series of an ample and globally generated line bundle $B$ is birational onto its image and $h^0(B)\geq7$ (see Theorem 1.11, \cite{GP1}). Recently Niu proved an analogue of Theorem \hyperref[SDR]{0.4} in dimension four. In fact he proved a general result for smooth, projective, regular, $K_X$-trivial varieties in all dimensions  (see \cite{Niu}) with an additional assumption of $H^2(\mathscr{O}_X)=0$ which in dimensions 2 and 3 recovers, and in dimension 4 generalizes the results of Saint-Donat and Gallego-Purnaprajna.\\
\indent We see that it is a natural question to ask whether and to what extent these theorems generalize to the other class of higher dimensional analogues of $K3$ surfaces, namely hyperk\"ahler varieties. Recall that there are many families of examples for Calabi-Yau varieties but only few classes of examples for hyperk\"ahler varieties are known upto deformation. Hilbert scheme of $n$ points on a $K3$ surface (we will denote them by $K3^{[n]})$, generalized Kummer varieties (we will denote them by $K^n(T)$) and two examples in dimension six and ten given by O'Grady  (we will denote it by $\mathscr{M}_6$ and $\mathscr{M}_{10}$ respectively) are the only known classes of examples upto deformation. Our main theorem for hyperk\"ahler varieties is the following.\par 

\noindent {\bf Theorem B.} \phantomsection\label{ob}\textit{Let $X$ be a projective hyperk\"ahler variety of dimension $2n\geq 4$ that is deformation equivalent to $K3^{[n]}$ or $K^n(T)$ or $\mathscr{M}_6$. Let $B$ be an ample and globally generated line bundle on $X$. Then the following happens.
\begin{itemize}
    \item[(1)] $B^{\otimes l}$ is projectively normal for $l\geq 2n$.
    \item[(2)] $B^{\otimes 2n-1}$ is projectively normal unless:
    \begin{itemize}
        \item [(a)] $n=2$, $X=K3^{[2]}$ and $\varphi_B$ maps $X$ onto a quadric (possibly singular) inside $\mathbb{P}^5$. In this case $q_X(B)=2$, $deg(\varphi_B)=6$, or
        \item[(b)] $n=3$, $X=K3^{[3]}$ and $\varphi_B$ maps $X$ onto a variety of minimal degree inside $\mathbb{P}^9$ which is obtained by taking cones over the Veronese embedding of $\mathbb{P}^2$ inside $\mathbb{P}^5$. In this case $q_X(B)=2$, $deg(\varphi_B)=30$.
    \end{itemize}
\end{itemize}
Hence if $X$ is as above and $B$ does not satisfy cases $2(a)$ or $2(b)$ then a general curve section of $|B|$ is non-hyperelliptic.}
\par 
As before, the study of the situation when $X$ maps onto a variety of minimal degree by the complete linear series $|B|$ is the main ingredient of our proof. We also use two key characteristics of a hyperk\"ahler variety which are the Riemann-Roch expression that comes from the existence of a primitive integral quadratic form on the second integral cohomology group of the variety and Matsushita's theorem on fibre space structure of a hyperk\"ahler manifold (see Theorem \hyperref[mat]{1.15}).\\
\indent Now we give the structure of our paper. In Section \hyperref[prelims]{2}, we briefly recall the main theorems and observations that we will use to prove our main results. Section \hyperref[CY]{3} and Section \hyperref[HK]{4} deals with the proof of Theorem \hyperref[oa]{A} and \hyperref[ob]{B} respectively.
\subsection*{Acknowledgements.} We want to thank our advisor Prof. B.P. Purnaprajna for his encouragement, support and guidance without which this work would have never been possible.

\section{Preliminaries and notations}\label{prelims}
Throughout this article, $X$ will always denote a smooth, projective variety  over $\mathbb{C}$. $K$ or $K_X$ will denote its canonical bundle. We will use the multiplicative and the additive notation of line bundles interchangeably. Thus, for a line bundle $L$, $L^{\otimes r}$ and $rL$ are the same. We have used the notation $L^{-r}$ for $(L^{*})^{\otimes r}$. We will use $L^{r}$ to denote the intersection product.

\subsection{Projectve normality.} For a globally generated line bundle $L$ on a smooth projective variety $X$, we have the following short exact sequence.
\begin{align}
    0 \to M_L \to H^0(L)\otimes \mathscr{O}_X \to L \to 0.\label{0}
\end{align}
First we recall the definition of projective normality and the $N_p$ property. 
\begin{definition}
Let $L$ be a very ample line bundle on a variety $X$. Let the following be the minimal graded free resolution of the coordinate ring $R$ of the embedding of $X$ induced by $L$
\begin{align*}
    0 \to F_n \to F_{n-1} \to \dots \to F_0\to R \to 0.
\end{align*}
Let $\mathscr{I}_X$ be the ideal sheaf of the embedding.
\begin{itemize}
    \item[(1)] $L$ satisfies the property $N_0$ (or embeds $X$ as a projectively normal variety) if $R$ is normal.
    \item[(2)] $L$ satisfies the property $N_1$ (or is normally presented) if in addition $\mathscr{I}_X$ is generated by quadrics.
    \item[(3)] $L$ satisfies the property $N_p$ if in addition to satisfying the property $N_1$, the resolution is linear from the second step until the $p$-th step.
\end{itemize}    
\end{definition}
We have the following necessary and sufficient condition for the $N_p$ property of an ample and base point free line bundle on $X$. In this article, we will only deal with the case when $p=0$.\par 
\phantomsection\label{critnp} 
\begin{theorem}
Let $L$ be an ample, globally generated line bundle on $X$. Suppose the cohomology group $H^1(\bigwedge^{p'+1}M_L\otimes L^{\otimes k})$ vanishes for all $0\leq p'\leq p$ and for all $k\geq 1$, then $L$ satisfies the property $N_p$. If in addition $H^1(L^{\otimes r})=0$ for all $r\geq 1$, then the above is a necessary and sufficient condition for $L$ to satisfy $N_p$.
\end{theorem}
We have made use of the  following observation of Gallego and Purnaprajna (see for instance \cite{GP1}) to show the surjections of certain multiplication maps.\par
\phantomsection\label{obs} 
\begin{observation}
Let $E$ and $L_1$, $L_2$, $\dots$, $L_r$ be coherent sheaves on a variety $X$. Consider the map  $H^0(E)\otimes H^0(L_1\otimes L_2\otimes\dots\otimes L_r)\xrightarrow[]{\psi}H^0(E\otimes L_1\otimes\dots\otimes L_r)$ and the following maps
$$H^0(E)\otimes H^0(L_1)\xrightarrow[]{\alpha_1}H^0(E\otimes L_1),$$
$$H^0(E\otimes L_1)\otimes H^0(L_2)\xrightarrow[]{\alpha_2}H^0(E\otimes L_1\otimes L_2),$$
$$\dots$$
$$ H^0(E\otimes L_1\otimes \cdots \otimes L_{r-1})\otimes H^0(L_r)\xrightarrow[]{\alpha_r}H^0(E\otimes L_1\otimes...\otimes L_r).$$
If $\alpha_1$, $\alpha_2$, $\dots$, $\alpha_r$ are surjective then $\psi$ is also surjective.
\end{observation}\par

The technique we use to show projective normalty of an ample and globally generated line bundle on a variety is  to use Koszul resolution to restrict the bundle on a smooth curve section and then showing the surjectivity of an appropriate multiplication map. It is worth mentioning that Koszul resolution is the special case of a particular complex, known as Skoda complex (see \cite{La}).\par 
\phantomsection\label{koszul}
\begin{definition}
Let $X$ be a smooth projective variety of dimension $n\geq 2$. Let $L$ be an ample and globally generated line bundle on $X$.
\begin{itemize}
    \item[(1)] Take $n-1$ general sections $s_1,\dots,s_{n-1}$ of $H^0(L)$ so the intersection of the divisor of zeroes $L_i=(s_i)_0$ is a nonsingular projective curve $C$, that is $C=L_1\cap \dots\cap L_{n-1}$.
    \item[(2)] Let $\mathscr{I}$ be the ideal sheaf of $C$ and let $W=span\{s_1,\dots,s_{n-1}\}\subseteq H^0(L)$ be the subspace spanned by $s_i$. Note that $W\subseteq H^0(L\otimes \mathscr{I})$.
The Koszul resolution ${\bf I}_1$ is the following exact sequence,
\begin{align*}
    0 \to \bigwedge\limits^{n-1}W\otimes L^{-(n-1)} \to \dots\to \bigwedge\limits^{2}W\otimes L^{-2} \to W\otimes L^{-1}\to \mathscr{I} \to 0.
\end{align*}
\end{itemize}
\end{definition}
Once we boil down our problem to a problem on curve, we use the following two results. The first one is a result of Green (see \cite{G2}).  
\phantomsection\label{green} 
\begin{lemma}
Let $C$ be a smooth, irreducible curve. Let $L$ and $M$ be line bundles on $C$. Let $W$ be a base point free linear subsystem of $H^0(C,L)$. Then the multiplication map $W\otimes H^0(M)\rightarrow H^0(L\otimes M)$ is surjective if $h^1(M\otimes L^{-1})\leq dim(W)-2$.
\end{lemma}
The second one is known as Castelnuovo-Mumford lemma (see \cite{Mu}). It is one of the main tools that has been used extensively in studying the syzygies of a variety embedded in a projective space.\par
\phantomsection\label{cm} 
\begin{lemma}
Let $L$ be a base point free line bundle on a variety $X$ and let $\mathscr{F}$ be a coherent sheaf on $X$. If $H^i(\mathscr{F}\otimes L^{-i})=0$ for all $i\geq 1$ then the multiplication map $H^0(\mathscr{F}\otimes L^{\otimes i})\otimes H^0(L)\rightarrow H^0(\mathscr{F}\otimes L^{\otimes i+1})$ surjects for all $i\geq 0$.
\end{lemma}\par
One can obtain effective projective normality results using the above theorems and lemmas. In order to improve the bounds one can use the following version of a result of Green (see \cite{G1}, p.1089, (3)).
\begin{theorem}\phantomsection\label{mg}
Let $X$ be a regular projective variety of dimension $n$ for which the canonical bundle $K$ is ample and base point free. Moreover, assume $h^0(K)\geq n+2$. Let $\varphi_K$ be the morphism induced by the complete linear series $|K|$. If $\varphi_K(X)$ is not a variety of minimal degree then the multiplication map $H^0(K)\otimes H^0(K^{\otimes l})\to H^0(K^{\otimes l+1})$ surjects for $l\geq n$.
\end{theorem}
The following corollary is the precise version of what we will use in the subsequent sections. 
\begin{corollary}\phantomsection\label{cor}
Let $X$ be a regular variety of dimension $n\geq 3$ with trivial canonical bundle. Let $L$ be an ample and globally generated line bundle on $X$ and $\varphi_L$ be the morphism induced by the complete linear series $|L|$. Assume $\varphi_L$ does not map $X$ onto a variety of minimal degree. Then $L^{\otimes n-1}$ is projectively normal.
\end{corollary}
\noindent\textit{Proof.} We use Theorem \hyperref[critnp]{1.2} and Observation \hyperref[obs]{1.3} to notice that to prove the projective normality of $L^{\otimes n-1}$, it is enough to show the surjectivity of the map $H^0(lL) \otimes H^0(L) \longrightarrow H^0((l+1)L)$ for $l\geq n-1$. We just prove the case for $l=n-1$ here, the rest follow similarly.\\
\indent Choose a smooth threefold section $T$ of the ample and base point free line bundle $L$. We have the following commutative diagram, \[
\begin{tikzcd}[row sep=large, column sep=2ex]
 0\arrow[r] & H^0((n-1)L)\otimes H^0(\mathscr{O}_X) \arrow[r] \arrow[d] &  H^0((n-1)L)\otimes H^0(L) \arrow[r] \arrow[d] &  H^0((n-1)L)\otimes H^0(K_T) \arrow[r]\arrow[d] & 0\\
0 \arrow[r] & H^0((n-1)L) \arrow[r] & H^0(nL) \arrow[r] & H^0(nK_T) \arrow[r] & 0
\end{tikzcd}\]
where $K_T = B|_T$ (by adjunction) denotes the canonical bundle of $T$. Since the leftmost map surjects, the middle vertical map surjects if and only if rightmost vertical map surjects. By Kodaira vanishing we have that $H^1((n-2)L)=0$ and hence we have a surjection,
\begin{equation*}
    H^0((n-1)L)\to H^0((n-1)L|_T)=H^0((n-1)K_T).
\end{equation*} 
Hence the rightmost vertical map surjects if and only if we have the surjection of $$H^0((n-1)K_T)\otimes H^0(K_T) \longrightarrow H^0(nK_T).$$ 
We notice that $T$ is a smooth irreducble regular variety of general type and hence by Theorem \hyperref[green]{1.7}, the above map  surjects unless  $T$ is mapped to a variety of minimal degree by the complete linear series of $L|_T$. But the latter is equivalent to saying that $X$ is mapped by the complete linear series of $|L|$ to a variety of minimal degree.\QEDB\par 
We briefly recall some basic definitions and results concerning varieties of minimal degree.
\subsection{Varieties of minimal degree} For any non-degenerate variety $X\subset \mathbb{P}^r$, we have the inequality 
\begin{equation*}
    deg(X)\geq 1+codim(X).
\end{equation*}
A variety $X\subset \mathbb{P}^r$ is said to be a \textit{variety of minimal degree} if it satisfies the equality $deg(X)=1+codim(X)$. If $codim(X)=1$ then of course the variety $X$ is a quadric hypersurface. The following theorem of Eisenbud and Harris (see \cite{Ei}) provides a complete classification of the varieties of minimal degree.
\begin{theorem}\phantomsection\label{classvmd}
Let $X\subset\mathbb{P}^r$ is a variety of minimal degree. Then $X$ is a cone over a smooth such variety. Moreover, if $X$ is smooth and $codim(X)>1$ then $X\subset\mathbb{P}^r$ is either a rational normal scroll or a Veronese surface $\mathbb{P}^2\subset\mathbb{P}^5$.
\end{theorem}
Recall that a rational normal scroll $X\subset\mathbb{P}^r$ of dimension $n$ is the image of a projective bundle $\pi:\mathbb{P}(\mathscr{E})\to\mathbb{P}^1$ through the morphism goven by the tautological line bundle $\mathscr{O}_{\mathbb{P}(\mathscr{E})}(1)$ where the vector bundle $\mathscr{E}=\mathscr{O}(a_1)\oplus\dots\oplus\mathscr{O}(a_n)$ satisfies $0\leq a_1\leq a_2\leq \dots\leq a_n$ and $deg(X)=\sum a_i$.\\
\indent If $a_1=a_2=\dots=a_l=0$ for some $0<l<n$ then $X$ singular and it is a cone over a smooth rational normal scroll. The vertex or singular locus $V$ of this cone has dimension $l-1$ and let $X_S=X\backslash V$ be the smooth part of $X$. Moreover, $X$ is normal and $\Tilde{X}=\mathbb{P}(\mathscr{E})\to X$ is a rational resolution of singularity which is called the $canonical$ $resolution$ of the rational normal scroll $X$.\\
\indent Let $X$ be a rational normal scroll. Let $H$ be the class of a hyperplane section and $R$ be the class of a general linear subspace of codimension one. We note the following fact: 
\begin{lemma}\phantomsection\label{rita}
If $codim(V,X)=2$ then $H\sim deg(X)R$.
\end{lemma}
\noindent\textit{Proof.} See \cite{Her}, Corollary 2.2 (2).

We prove a lemma here that we will use in the subsequent sections to prove our main theorems.
\begin{lemma}\phantomsection\label{codimv}
Let $X$ be a smooth projective variety and let $L$ be an ample and globally generated line bundle on $X$. Let $\varphi_L$ be the morphism induced by the complete linear series $L$. Assume $\varphi_L$ maps $X$ onto a singular rational normal scroll $Y$ with vertex $V$. Then $codim(V,Y)=2$ that is $Y$ is obtained by taking cones over a rational normal curve.
\end{lemma}
\noindent \textit{Proof.} Suppose $\Tilde{Y}=\mathbb{P}(\mathscr{E})\to Y$ is the canonical resolution. If the codimension of the singular locus is $> 2$, then the canonical resolution is a small resolution (see \cite{Her}, Proposition 2.1) and hence $X$ could be obtained by contracting subschemes of $(X \times_Y \Tilde{Y})_{red}$ of codimension greater than one, which contradicts the factoriality of $X$. The assertion follows since $codim(V,Y)\neq 1$ as $Y$ is normal.\QEDB

\subsection{Hyperk\"ahler varieties} Recall that a hyperk\"ahler variety $X$ is a compact K\"ahler manifold that is simply connected, projective and its space of global holomorphic two forms is spanned by a symplectic form. The symplectic form ensures that $K_X$ is trivial and $dim(X)$ is even. It is also known that $\chi(\mathscr{O}_X)=n+1$ and the following are the values of $h^p(X,\mathscr{O}_X)$,
\[h^p(X,\mathscr{O}_X)=
    \begin{dcases}
        1 & \textrm{$p$ is even,} \\
        0 & \textrm{$p$ is odd.} \\
    \end{dcases}
\]
\indent Only a few classes of examples of hyperk\"ahler varieties are known. Beauville first gave examples of two distinct deformation classes of compact hyperk\"ahler manifolds in all even dimensions greater than or equal to $2$ (see \cite{Be}). The first example is the Hilbert scheme $S^{[n]}$ of length $n$ subschemes on a $K3$ surface $S$. The second one is the generalized Kummer variety $K^n(T)$ which is the fibre over the $0$ of an Abelian variety $T$ under the morphism $\phi\circ\psi$ (see the diagram below)
\begin{align*}
    T^{[n+1]} \xrightarrow{\psi} T^{(n+1)}\xrightarrow{\phi} T
\end{align*}
where $T^{[n+1]}$ Hilbert scheme of length $n+1$ subschemes on the Abelian variety $T$, $T^{(n+1)}$ is the symmetric product, $\psi$ is the Hilbert chow morphism and $\phi$ is the addition on $T$. Two other distinct deformation classes of hyperkahler manifolds $\mathscr{M}_6$ and $\mathscr{M}_{10}$ are given by O'Grady in dimensions $6$ and $10$ respectively which appear as desingularizations of certain modulii spaces of sheaves over symplectic surfaces (see \cite{OG1}, \cite{OG2}). All other known examples are deformation equivalent to one of these.\\
\indent We start by the following theorem of Beauville and Fujiki (see \cite{Be} and \cite{Fu}).\par
\phantomsection\label{cq}
\begin{theorem}
Let $X$ be a hyperk\"ahler variety of dimension $2n$. There exists a quadratic form $q_X:H^2(X,\mathbb{C})\rightarrow \mathbb{C}$ and a positive constant $c_X\in\mathbb{Q}_+$ such that for all $\alpha$ in $H^2(X,\mathbb{C})$, $\int_X\alpha^{2n}=c_X\cdot q_X(\alpha)^n$. The above equation determines $c_X$ and $q_X$ uniquely if one assumes the following two conditions.
\begin{itemize}
    \item[(I)] $q_X$ is a primitive integral quadratic form on $H^2(X,\mathbb{Z})$;
    \item[(II)] $q_X(\sigma,\bar{\sigma})>0$ for all $0\neq\sigma\in H^{2,0}(X)$.
\end{itemize}
Here $q_X$ and $c_X$ are called the Beauville form and Fujiki constant respectively.
\end{theorem}
The Beauville form and Fujiki constants are fundamental invariants of a hyperk\"ahler variety. They play an important role in determining the intersections on $X$ as the following theorem shows (See \cite{Fu}, \cite{Gr}).
\phantomsection\label{defc}
\begin{theorem}
Let $X$ be a hyperk\"ahler variety of dimension $2n$. Assume that $\alpha\in H^{4j}(X,\mathbb{C})$ of type $(2j,2j)$ on all small deformations of $X$. Then there exists a constant $C(\alpha)\in \mathbb{C}$ depending on $\alpha$ such that $\int_X\alpha\cdot\beta^{2n-2j}=C(\alpha)\cdot q_X(\beta)^{n-j}$ for all $\beta\in H^2(X,\mathbb{C})$.
\end{theorem}
\phantomsection\label{grr}As a consequence of the theorem above, we get the following form of the Riemann-Roch formula for an line bundle $L$ on a hyperk\"ahler variety of dimenson $2n$ (see \cite{Hu}), 
\begin{align}
    \chi(X,L)=\sum_{i=0}^{n}\dfrac{a_i}{(2i)!}q_X(c_1(L))^i
\end{align} where $a_i=C(td_{2n-2i}(X))$. Here $a_i$'s are constants depending only on the topology of $X$.\\
\indent Elingsrad-Gottsche-Lehn computes the rational constants of the Riemann roch expression for hyperk\"ahler manifolds of deformation type $K3^{[n]}$ (See \cite{Ell}) and Britze-Nieper computes the same for generalized Kummer varieties $K^n(T)$ of dimension $2n$ (see \cite{Br}). If $X$ is of $K3^{[n]}$ type we have,
\begin{align}
    \chi(L) = \binom{\frac{1}{2}q(L)+n+1}{n}\quad and\quad c_X=\dfrac{(2n)!}{n!2^n}.\label{rrk3n}
\end{align}
For a generalized Kummer variety of dimension $2n$ we have the following Riemann-Roch formula,
\begin{align}
    \chi(L) = (n+1) \binom{\frac{1}{2}q(L)+n}{n}\quad and \quad c_X=(n+1)\dfrac{(2n)!}{n!2^n}\label{rrknt}
\end{align}
The Riemann-Roch formula and Fujiki constant for $\mathscr{M}_6$ are the same as that of $K^3(T)$.\\
\indent We will use Matsushita's theorem on fibre space structure of hyperk\"ahler varieties. We recall the the definition and the main theorem.
\begin{definition}
Let $X$ be an algebraic variety. A fibre space structure of $X$ is a proper surjective morphism $f:X\to S$ that satisfies the following two conditions:
\begin{itemize}
    \item[(1)] $X$ and $S$ are normal varieties with $0<dim(S)<dim(X)$.
    \item[(2)] A general fibre of $f$ is connected.
\end{itemize}
\end{definition}\phantomsection\label{mat}
\indent The result that we will use in Section \hyperref[HK]{4} is the following (see \cite{Ma}, Theorem 2, (3)).
\begin{theorem}
Let $f: X\to B$ is a fibre space structure on a projective hyperk\"ahler variety $X$ of dimension 2n with projective base $B$. Then $dim(B)=n$.
\end{theorem}

\section{Proof of Theorem A}\label{CY}
 The main aim of this section is to prove results on effective very ampleness and projective normality on a four dimensional variety with trivial canonical bundle. We start with a general statement on projective normality and normal presentation.\par 

\begin{lemma}
\phantomsection\label{old} Let $X$ be a smooth, projective $n$-fold with trivial canonical bundle. Let $B$ be an ample and base point free line bundle on $X$. Let $h^0(B)\geq n+2$. Then $lB$ satisfies the property $N_0$ for all $l\geq n$. Moreover, if $X$ is Calabi-Yau, then $lB$ satisfies the property $N_1$ for all $l\geq n$.
\end{lemma}

\noindent\textit{Proof.} Follows immediately from Theorem 2.3 and Theorem 3.4 of \cite{MR}. \QEDB\par 

Now we want to find out what multiple of an ample line bundle is very ample on a four dimensional variety with trivial canonical bundle. We will use the Fujita freeness on four folds that has been proved by Kawamata in \cite{Ka}. We begin with a lemma.\par

\begin{lemma}
\phantomsection\label{surjmultcy}Let $X$ be a fourfold with trivial canonical bundle. Let $A$ be an ample line bundle and let $B=nA$ for $n\geq 5$. Then the multiplication map
$H^0(3B+kA)\otimes H^0(B) \rightarrow H^0(4B+kA)$ is surjective for $k\geq1$.
\end{lemma} 

\noindent\textit{Proof.} Note that $B$ is base point free by Kawamata's proof of Fujita's base point freeness theorem on fourfolds (see \cite{Ka}). We prove the statement for $k=1$. For $k>1$ the proof is exactly the same.\\
\indent Let $C$ be a smooth and irreducible curve section of the linear system $|B|$ and let $\mathscr{I}$ be the ideal sheaf of $C$ in $X$.
We have the following commutative diagram with the two horizontal rows exact. Here $\mathscr{I}$ is the ideal sheaf of $C$ in $X$, $V$ is the cokernel of the map $H^0(B\otimes \mathscr{I}) \rightarrow H^0(B)$. \[
\begin{tikzcd}[row sep=large, column sep=2ex]
    0 \arrow{r} & H^0(B\otimes \mathscr{I})\otimes H^0(3B+A) \arrow{r} \arrow{d} & H^0(B)\otimes H^0(3B+A) \arrow{r} \arrow{d} & V\otimes H^0(3B+A) \arrow{r} \arrow{d} & 0  \\
    0 \arrow{r} & H^0((4B+A)\otimes \mathscr{I}) \arrow{r}  & H^0((4B+A)) \arrow{r}  &  H^0(4B+A|_C) \arrow{r}  & 0  
\end{tikzcd} \]
\indent Now we claim that the leftmost vertical map is surjective. Consider the Koszul resolution,
\begin{align}
    0 \to \bigwedge\limits^{3} W\otimes B^{-3} \to \bigwedge\limits^2 W\otimes B^{-2} \to W\otimes B^* \to \mathscr{I} \to 0.
\end{align}
Tensor it with $4B+A$ to get the following,
\begin{align}
    0 \to \bigwedge\limits^3 W\otimes (B+A) \xrightarrow{f_3} \bigwedge\limits^2 W\otimes (2B+A) \xrightarrow{f_2} W\otimes (3B+A) \xrightarrow{f_1} (4B+A)\otimes \mathscr{I} \to 0  
\end{align}
That gives us two short exact sequences.
\begin{align}
    0 \to Ker(f_1) \to W\otimes (3B+A) \xrightarrow{f_1} (4B+A)\otimes \mathscr{I} \to 0,
\end{align}
\begin{align}
    0 \to\bigwedge\limits^3 W\otimes (B+A) \xrightarrow{f_3} \bigwedge\limits^2 W\otimes (2B+A) \xrightarrow{f_2} Ker(f_1) \to 0.
\end{align}
Taking long exact sequence of cohomology in the second sequence we get the following,
\begin{align}
    \bigwedge\limits^2 W\otimes H^1(2B+A) \to H^1(Ker(f_1)) \to  \bigwedge\limits^3 W\otimes H^2(B+A).
\end{align}
Hence $H^1(Ker(f_1))=0$ since the other terms of the exact sequence vanish by Kodaira Vanishing.\\
The long exact sequence of cohomology associated to the first sequence is the following,
\begin{align}
    W\otimes H^0(3B+A) \to H^0((4B+A)\otimes \mathscr{I})) \to H^1(Ker(f_1)).
\end{align}
We showed that the last term is zero, thus $W\otimes H^0(3B+A) \rightarrow  H^0((4B+A)\otimes \mathscr{I}))$ surjects. Consequently, $H^0(B\otimes \mathscr{I})\otimes H^0(3B+A) \rightarrow H^0((4B+A)\otimes \mathscr{I})$ surjects since $W \subseteq H^0(B\otimes \mathscr{I})$.\\
\indent In order to prove the lemma we are left to show that $V\otimes H^0(3B+A) \rightarrow H^0((4B+A)|_C)$ surjects.
Since we have the surjection of $H^0(3B+A) \rightarrow H^0((3B+A)|_C) $,
 it is enough to show the surjection of $V\otimes H^0(3B+A|_C) \rightarrow H^0((4B+A)|_C)$. Thus, using Lemma \hyperref[green]{1.5} we need to prove that 
 \begin{align*}
     h^1(2B+A|_C)\leq dim V-2.
 \end{align*}
\indent To prove this inequality, first we tensor the Koszul resolution by $B$,
\begin{align}
    0 \to \bigwedge\limits^3 W\otimes B^{-2} \xrightarrow{f_3}  \bigwedge\limits^2 W\otimes B^* \xrightarrow{f_2} W\otimes \mathscr{O}_X \xrightarrow{f_!} B\otimes \mathscr{I} \to 0.
\end{align}
As before, we end up getting two short exact sequences,
\begin{align}
    0 \to Ker(f_1) \to W\otimes \mathscr{O}_X \xrightarrow{f_1} B\otimes \mathscr{I} \to 0 
\end{align}
\begin{align}
    0 \to \bigwedge\limits^3 W\otimes B^{-2} \xrightarrow{f_3} \bigwedge\limits^2 W\otimes B^* \xrightarrow{f_2} Ker(f_1) \to 0
\end{align}
The long exact sequence of cohomology associated to the second sequence gives,
\begin{align}
    \bigwedge\limits^2 W\otimes H^0(B^*) \to H^0(Ker(f_1)) \to \bigwedge\limits^3 W\otimes H^1(B^{-2}).
\end{align}
Consequently, $H^0(Ker(f_1))=0$ since $H^1(B^{-2})=0$ by Kodaira vanishing and $H^0(B^*)=0$. \\Taking cohomology once more we have the following exact sequence,
\begin{align}
    \bigwedge\limits^2 W\otimes H^1(B^*) \to H^1(Ker(f_1)) \to  \bigwedge\limits^3 W\otimes H^2(B^{-2}).
\end{align}
Hence $H^1(Ker(f_1))=0$ since the other terms of the exact sequence vanish by Kodaira Vanishing. \\
The long exact sequence of cohomology associated to the first sequence is the following.
\begin{align}
     H^0(Ker(f_1)) \to W\otimes H^0(\mathscr{O}_X) \to H^0(B\otimes \mathscr{I})) \to H^1(Ker(f_1))
\end{align}
But the first and last terms are zero by Kodaira Vanishing and hence $h^0(B\otimes \mathscr{I})=dim W\leq 3$. Thus we obtain the inequality $dim V -2\geq h^0(B)-5$.\\
\indent On the other hand the canonical bundle of $C$ is given by $3B|_C$. Applying Serre-Duality it is enough to prove that $ h^0(B-A)\leq h^0(B)-5$
i.e. $ h^0((n-1)A)\leq h^0(nA)-5$.\\
\indent Applying Riemann Roch for $nA$ and $(n-1)A$ and subtracting the equations we obtain,
\begin{align}
    h^{0}(nA) - h^{0}((n-1)A) = \displaystyle{\frac{n^4-(n-1)^4}{24}}A^4+\displaystyle{\frac{n^2-(n-1)^2}{24}}A^2\cdot c_2.
\end{align}
The result of Miyaoka (see \cite{Mi}) shows that $A^2\cdot c_2 \geq 0$ which gives,
\begin{align}
    h^{0}(nA) - h^{0}((n-1)A) \geq \displaystyle{\frac{n^4-(n-1)^4}{24}}A^4 \geq 5 \textrm{ if $n \geq 5$}
\end{align}
and that concludes the proof.\QEDB\par
 
  Now we give a classification theorem in which we classify the varieties which come as an image  of a regular fourfold with trivial canonical bundle by an ample, globally generated line bundle with an additional property of being a variety of minimal degree. \par
\begin{theorem}
\phantomsection\label{class}Let $X$ be a regular four-fold with trivial canonical bundle. Let $\varphi$ be the morphism induced by the complete linear series of an ample and base point free line bundle $B$ on $X$ with $h^0(B)=r+1$ and let $d$ be the degree of $\varphi$. If $\varphi$ maps $X$ to a variety of minimal degree $Y$ then,
\begin{align*}
    d \leq \displaystyle{\frac{24(r-1)}{r-3}}.
\end{align*}
(a) Assume $Y$ is smooth. Then one of the following happens:
\begin{itemize}
    \item[(1)] $Y=\mathbb{P}^4$.
    \item[(2)] $Y$ is a smooth quadric hypersurface in $\mathbb{P}^5$.
    \item[(3)] $Y$ is a smooth rational normal scroll of dimension $4$ in $\mathbb{P}^6$ or $\mathbb{P}^7$ and $X$ is fibered over $\mathbb{P}^1$. Moreover, the general fibre is a smooth threefold $G$ with $K_G=0$ and the degree $d$ satisfies the following bounds;
    $$2\leq d \leq min \bigg\{6h^0(B|_G), \displaystyle{\frac{24(r-1)}{r-3}}\bigg\}.$$
    If in addition $G$ is regular we have the following;
    \begin{center}
    $2h^0(B|_G)-6\leq d \leq min \bigg \{ 6(h^0(B|_G)-1), \displaystyle\frac{24(r-1)}{r-3}\bigg\}$, if $d$ is even and \\
    $2h^0(B|_G)-5\leq d \leq min \bigg\{6(h^0(B|_G)-1),\displaystyle\frac{24(r-1)}{r-3}$\bigg\}, if $d$ is odd.  \end{center}
    \item[(4)] $Y$ is a smooth rational normal scroll in $\mathbb{P}^r$ for $r\geq8$ and $X$ is fibered over $\mathbb{P}^1$ and the general fibre is a three-fold $G$ with $K_G=0$ and the degree $d$ of $\varphi_B$ satisfies $2\leq d \leq 18$.
\end{itemize}
(b) Assume $Y$ is singular. Then one of the following happens:
\begin{itemize}
    \item[(1)] $Y$ is a singular quadric hypersurface.
    \item[(2)] $Y$ is a singular four-fold which is either a triple cone over a rational normal curve in $\mathbb{P}^r$ where $6 \leq r \leq 8$ or a double cone over the Veronese surface in $\mathbb{P}^5$.
\end{itemize}
\end{theorem}

\noindent\textit{Proof.} We first prove the inequality. Using Riemann-Roch we can see that,
\begin{align}
    h^0(B)=\displaystyle{\frac{1}{24}}B^4+\displaystyle{\frac{1}{24}}B^2\cdot c_2+2.
\end{align}
and we also have that $B^4=d(r-3)$ since $Y$ is a variety of minimal degree. By Miyaoka's result (see \cite{Mi}) we have that $B^2\cdot c_2 \geq 0$ and hence we have the inequality $d \leq \frac{24(r-1)}{r-3}$.\par 

\noindent (a) We now describe the cases when $Y$ is a smooth variety of minimal degree. We have that $r\geq 4$. \par 

\noindent \underline{\textit{Case 1.}} If $r=4$, we have that $Y=\mathbb{P}^4$. \par
\noindent \underline{\textit{Case 2.}} If $r=5$, we have that codimension of $Y$ is one and degree is $2$ which implies that $Y$ is a smooth quadric hypersurface.\par

 Suppose $r\geq 6$, we have that $Y$ is a smooth rational normal scroll (which is abstractly a projective bundle over $\mathbb{P}^1$) and is hence fibered over $\mathbb{P}^1$. Let this map from $Y$ to $\mathbb{P}^1$ be $\phi$. Composing this with $\varphi$ we get a map $\phi \circ \varphi$ : $X \rightarrow \mathbb{P}^1$. Therefore, $X$ is fibered over $\mathbb{P}^1$. The general fibre is the inverse image of the general linear fiber of the smooth scroll  and is hence irreducible by Bertini's theorem (see \cite{FL}, Theorem 1.1). This along with generic smoothness implies that the general fibre of $\phi$ is a smooth threefold $G$ with $K_G=0$ by adjunction. Let the general fibre of $Y$ be denoted by $R$ and that of $X$ is denoted by $G$. We have the following exact sequence of cohomology of line bundles on $X$.
\begin{align}
    0 \to H^0(B(-G)) \to H^0(B) \to H^0(B\otimes \mathscr{O}_G) \to H^1(B(-G)). 
\end{align}
Notice that $H-R$ is a base point free divisor in $Y$ where $H$ is a hyperplane section in $Y$. We have that $Y$ is $S(a_1,a_2,a_3,a_4)$ i.e, $Y$ is the image of $\mathbb{P}(\mathscr{E})$ where $\mathscr{E}$ is the following vector bundle,
\begin{align*}
    \mathbb{P}(\mathscr{E})=\mathscr{O}(a_1)\oplus \mathscr{O}(a_2)\oplus \mathscr{O}(a_3) \oplus \mathscr{O}(a_4)
\end{align*}
mapped to the projective space by $|\mathscr{O}_{\mathbb{P}(\mathscr{E})}(1)|$.\par 
\noindent\underline{\textit{Case 3.}} For the cases $r=6$ or $r=7$ we use the fact that degree of $\varphi$ is equal to the degree of $\varphi|_G$ and then use Riemann-Roch theorem on the threefold $G$ (see \cite{Har}, Appendix A, Exercise 6.7) noticing the fact that $K_G=0$ and that $B|_G$ is ample and base point free. This gives the upper bound $6h^0(B|_G)$ since we have that $B|_G\cdot c_2 \geq 0$ (see \cite{Mi}). The lower bound $2$ is due to the fact that $G$ cannot be birational to $\mathbb{P}^3$.\\
\indent Assuming $G$ is regular and hence Calabi-Yau we have that $h^0(B|_G) \geq \frac{1}{6}(B|_G)^3+1$ and hence we have $d \leq 6(h^0(B|_G)-1)$. The lower bound is obtained by Proposition 2.2, part (1) of \cite{Kan}.\par 
\noindent\underline{\textit{Case 4.}} Suppose $r\geq 8$. Recall that $H-R$ is base point free.  We compute $(H-R)^4=H^4-4H^3R$;
\begin{align}
    H^4=\sum_{i=0}^{3}a_iH^3R \textrm{ and }r=\sum_{i=0}^{3}a_i+3.
\end{align}
So, $r\geq 8$ gives $\sum_{i=0}^{3}a_i \geq 5$ which gives $(H-R)^4>0$ as $H$ is ample. Hence, $H-R$ is nef and big and consequently $B-G$ is nef and big as well. Thus by Kawamata-Viehwag vanishing, we have that $H^1(B-G)=0$. Hence $\varphi|_G$ is given by the complete linear system $|B|_G|$. \\
\indent Since $G$ maps to $F=\mathbb{P}^3$ we have that $h^0(B|_G)=4$. Now, the degree of $\varphi$ is also the degree of $\varphi|_G$ for a general fibre $G$. Hence by a result of Gallego and Purnaprajna (see \cite{GP1}, Theorem 1.6) we have that $2 \leq d \leq 18$. \par
\noindent(b) Now we assume that $Y$ is singular.\par 
\noindent\underline{\textit{Case 1.}} As before, if $r=5$ then $Y$ is a singular quadric.

\noindent\underline{\textit{Case 2.}} Suppose the image $Y$ of $X$ under the morphism defined by $|B|$ is a singular variety. Then by Theorem \hyperref[classvmd]{1.8}, $Y$ is a cone over a smooth variety of minimal degree with vertex $V$. Moreover, $codim(V,Y)=2$ by Lemma \hyperref[codimv]{1.11}. Hence $Y$ can be either a triple cone over a rational normal curve or a double cone over the Veronese surface in $\mathbb{P}^5$.\\
\indent Suppose $Y$ is a cone over a rational normal curve. Let $R$ be a general linear subspace of $Y$ and $G$ its inverse image. By Bertini (see \cite{FL}, Theorem 1.1), $G$ is irreducible and $d = B^3\cdot G$. Moreover, using the fact that the codimension of the singular locus of $Y$ is exactly $2$, we have $B = deg(Y) \cdot G = (r-3)G$ by Lemma \hyperref[rita]{1.9}. Hence $d \geq (r-3)^3$. Using the previously proved upper bound we have, $$ (r-3)^3 \leq 24(r-1)$$ and hence $r \leq 8$.\QEDB
 
 Now we prove Theorem \hyperref[oa]{A} using the previous theorem. We notice that part $(2)$ of the next theorem holds for both hyperk\"ahler and Calabi-Yau four-folds in dimension four since it only requires a regular fourfold with trivial canonical bundle.\par
 
\begin{theorem}
\phantomsection\label{pthma}Let $X$ be a four-fold with trivial canonical bundle and let $A$ be an ample line bundle on $X$. then, 
\begin{itemize}
    \item[(1)] $nA$ is very ample and embeds $X$ as a projectively normal variety for $n \geq 16$.
    \item[(2)] If $H^{1}(\mathscr{O}_{X})=0$ then $nA$ is very ample and embeds $X$ as a projectively normal variety for $n \geq 15$.
\end{itemize}
\end{theorem} 
\noindent\textit{Proof of (1).} By the result of Kawamata (see \cite{Ka}), we have that on a fourfold with trivial canonical bundle if $A$ is ample then $nA$ is base point free for $n \geq 5$. Now using CM lemma (see Lemma \hyperref[cm]{ 1.6}) we can easily prove that $nA$ satisfies the property $N_0$ for $n \geq 21$.\par 
\indent If we set $B = 5A$ then $20A = 4B$ and it satisfies the property $N_0$ by Lemma \hyperref[old]{2.1}.\\ 
\indent Using Lemma \hyperref[surjmultcy]{2.2}, Lemma \hyperref[cm]{1.6} and Observation \hyperref[obs]{1.3}, we can see that the following map
\begin{align*}
    H^{0}(nkA) \otimes H^{0}(nA) \longrightarrow H^{0}((nk+n)A)
\end{align*}
is surjective for $k \geq 2$ and $16 \leq n \leq 19$. 
 So we are left to check the surjectivity of  the multiplication map
 $H^{0}(nA) \otimes H^{0}(nA) \longrightarrow H^{0}(2nA)$ for $16 \leq n \leq 19$. We just prove it for $n=16$. The proof is similar for the other three cases.\\
\indent For $n = 16$, the following maps surjects by Lemma \hyperref[surjmultcy]{2.2} and Lemma 
\hyperref[cm]{1.6},
 \begin{align}
     H^{0}(16A) \otimes H^{0}(5A) \to H^{0}(21A)\textrm{ and }H^{0}(21A) \otimes H^{0}(5A) \to H^{0}(26A).
 \end{align}
 \noindent Therefore, by Observation \hyperref[obs]{1.3} we need to show that  $H^{0}(26A) \otimes H^{0}(6A) \longrightarrow H^{0}(32A)$ is surjective which follows from Lemma \hyperref[cm]{1.6} as well. \par
 
 \noindent\textit{Proof of (2).} Suppose $H^1(\mathscr{O}_X)=0$. We just need to show that $15A$ satisfies the property $N_0$.\\
 \indent Let $B=5A$ which is ample and base point free (see \cite{Ka}). By Corollary \hyperref[cor]{1.8}, we know that $3B$ is projectively normal unless the image of the morphism induced by the complete linear series $|B|$ is a variety of minimal degree. Thus, we aim to show that the image of the morphism induced by the complete linear series $|5A|$ is not a variety of minimal degree. Applying Riemann-Roch we get,
 \begin{align}
     h^{0}(5A) = \displaystyle{\frac{625}{24}}A^4 + \displaystyle{\frac{25}{24}}A^2\cdot c_2+2 \geq 28.
 \end{align}
Now suppose that the image is a variety of minimal degree. However, since the codimension of the image is $\geq 24$, by Theorem \hyperref[class]{2.3}, we have that the image cannot be a quadric hypersurface or a cone over the veronese embedding of $\mathbb{P}^2$ in $\mathbb{P}^5$ or a cone over a rational normal curve.
  Hence the image is a smooth rational normal scroll.\\
\indent  Let $h^0(B)=r+1$. Hence the degree of the image is $r-3$. 
  Also, let the degree of the finite morphism given by the complete linear series of $B$ be $d$. We know by Theorem  \hyperref[class]{2.3} that 
  \begin{align*}
      d\leq \frac{24(r-1)}{r-3}.
  \end{align*}
  Using $h^0(B)\geq 28$ we have that $r\geq 27$ and hence $d\leq 26$.\\
  \indent Since the image of the morphism is a smooth rational normal scroll of dimension $4$, we can choose a general $\mathbb{P}^3=R$ and take the pullback of the divisor $R$ under the morphism induced by the complete linear series $|B|$ and call it $G$. The degree of the morphism restricted to $G$ is again $d$.
  Since the degree of $R$ in the image is $1$ we have that $d=B^3\cdot G=125A^2\cdot G\geq 125$ (since $A$ is ample and $G$ is effective) contradicting $d\leq 26$. Hence the image cannot be variety of minimal degree. \QEDB\par

\section{Proof of Theorem B}\label{HK}
 
Let $X$ be a hyperk\"ahler variety of dimension $2n$ and let $B$ be an ample and globally generated line bundle on $X$. In this section, $\varphi_{B}$ will always denote the morphism induced by the complete linear series $|B|$. The aim is to study the projective normality of $B^{\otimes 2n-1}$. We do this by Theorem \hyperref[cor]{1.8} that require us to analyze the case when $\varphi_B$ maps $X$ onto a variety of minimal degree. Recall that a variety of minimal degree is either (1) a quadric hypersurface, or (2) a smooth rational normal scroll, or (3) cone over a smooth rational normal scroll, or (4) cone over the Veronese embedding of $\mathbb{P}^2$ inside $\mathbb{P}^5$ (\cite{Ei}). The following lemma eliminates a few cases.\par 
\begin{lemma}\phantomsection\label{elimination}
Let $B$ be an ample and globally generated line bundle on a hyperk\"ahler variety $X$ of dimension $2n$. Suppose the morphism $\varphi_B$ induced by the complete linear series $|B|$ maps $X$ onto a variety of minimal degree $Y$. Then,
\begin{itemize}
    \item[(1)] If $Y$ is a quadric hypersurface then $h^0(B)=2n+2$.
    \item[(2)] $Y$ can never be a smooth rational normal scroll.
    \item[(3)] If $Y$ is a cone over a smooth rational normal scroll then the codimension of its singular locus is two i.e. $Y$ is obtained by taking cones over a smooth rational normal curve.
    \item[(4)] If $Y$ is a cone over the Veronese embedding of $\mathbb{P}^2$ inside $\mathbb{P}^5$ then $h^0(B)=2n+4$.
\end{itemize}    
\end{lemma}
\noindent\textit{Proof.} (1) and (4) are obvious. (3) comes from Lemma \hyperref[codimv]{1.11}. We give the proof of (2) below.\\
\indent To prove (2) we argue by contradiction. Suppose the image is a smooth rational normal scroll. Since a smooth scroll admits a morphism to $\mathbb{P}^1$ we have a composed morphism from $X$ to $\mathbb{P}^1$. Take the Stein factorization of this morphism which has connected fibres and notice that since $X$ is smooth this further factors through a normalization. So we get a morphism from $X$ to a normal base of dimension $1$ (hence smooth in this case) with connected fibres which contradicts Matsushita's result on the fibre space structure of a holomorphic symplectic manifold (see Theorem \hyperref[mat]{1.15}).\QEDB\par 

We give two definitions below that we will use later in this note.  
\begin{definition}\label{RRR}
For a given hyperk\"ahler variety $X$, we define the following two polynomials,
\begin{equation*}
    RR_X(x)=\displaystyle\sum_{i=0}^n\dfrac{a_i}{(2i)!}(x)^i \quad and\quad R_X(x)=RR_X(x)-\dfrac{a_n}{(2n)!}x^n
\end{equation*}
where $a_i=C(td_{2n-2i})$. (Note: these polynomials depend only on the deformation type of $X$.)
\end{definition}
\begin{definition}
For a given hyperk\"ahler variety $X$ with Beauville form $q_X$, we define the constant $\alpha_X$ as below,
\begin{equation*}
    \alpha_X=min\bigg\{q_X(A)\bigg\vert \textrm{ A is an ample line bundle on $X$}\bigg\}.
\end{equation*}
\end{definition}
Our next task is to find an upper bound for $deg(\varphi_B)$ for an ample and base point free line bundle $B$. The next result in fact just uses the fact that $B$ is base point free and the induced morphism $\varphi_B$ is generically finite.
\begin{lemma}\phantomsection\label{ub}
Let $X$ be a hyperk\"ahler manifold of dimension $2n$. Assume $R_X\vert_{\mathbb{Z}\geq 0}$ is increasing and $R_X(\alpha_X)>2n$. Then, for any globally generated line bundle $B$ on $X$, such that $\varphi_B$ is generically finite, $deg(\varphi_B)<(2n)!$. 
\end{lemma}
\noindent\textit{Proof.} Let $Y=Im(\varphi_B)$. We have the equation $B^{2n}=deg(\varphi_B)\cdot deg(Y)$ to begin with. Note that,
\begin{equation}
    deg(Y)\geq 1+codim(Y)\implies deg(Y)\geq h^0(B)-2n
\end{equation}
since $codim(Y)=h^0(B)-2n$. Thus, we get,
\begin{equation}
    B^{2n}\geq deg(\varphi_B)(h^0(B)-2n).\label{equb}
\end{equation}
\indent Recall that $B^{2n}=c_X(q_X(B))^n$ and notice that $h^0(B)=RR_X(q_X(B))=R_X(q_X(B))+\dfrac{a_n}{(2n)!}(q_X(B))^n$. Since, $a_n=c_X$, using \hyperref[equb]{(3.2)} we get,
\begin{equation*}
    c_X(q_X(B))^n\geq deg(\varphi_B)\left(R_X(q_X(B))+\dfrac{a_n}{(2n)!}(q_X(B))^n-2n\right)
\end{equation*}
\begin{equation*}
    \implies c_X(q_X(B))^n\left(1-\dfrac{deg(\varphi_B)}{(2n)!}\right)\geq deg(\varphi_B)(R_X(q_X(B))-2n)\geq deg(\varphi_B)(R_X(\alpha_X)-2n).
\end{equation*}
That concludes the proof since the last term is strictly greater than zero by hypothesis and $q_X(B) > 0$ since $B$ is nef and big.\QEDB\par 
\begin{remark}
If all Todd classes of the hyperk\"ahler variety X is fakely effective then $R_X\vert_{\mathbb{Z}\geq 0}$ is increasing. In particular, it is satisfied for all known examples of hyperk\"ahler varieties, except O'Grady's 10 dimensional example $\mathscr{M}_{10}$ (see \cite{Cao}, Theorem 1.8) remaining unknown, which is also clear from their explicit Riemann-Roch expressions.\QEDB 
\end{remark}
The remark above leads to the following consequence.
\begin{remark}\phantomsection\label{hypcheck}
The hypothesis of Lemma \hyperref[ub]{3.4} are satisfied for all known examples of hyperk\"ahler varieties $X$ of dimension $2n\geq 4$ (see \cite{Cao}, Theorem 1.8), except O'Grady's 10 dimensional example $\mathscr{M}_{10}$.
\end{remark}
\noindent\textit{Proof.} Thanks to the previous remark, it is enough to show that $R_X(\alpha_X)>2n$.\\
\indent If $X$ is either type $K3^{[n]}$ or $K^n(T)$ then $\alpha_X\geq 2$. Since $R_X\vert_{\mathbb{Z}\geq 0}$ is increasing, we have, 
\begin{align}
    R_X(\alpha_X)\geq RR_X(2)-\dfrac{c_X2^n}{(2n)!}.
\end{align}
\indent If $X$ is of type $K3^{[n]}$ (resp. $K^n(T)$), using the Riemann-Roch expression \hyperref[rrk3n]{1.3} (resp. \hyperref[rrknt]{1.4}) and $n\geq 2$ we have the following,
\begin{align}
    RR_X(2)-\dfrac{c_X2^n}{(2n)!}=\binom{n+2}{n}-\dfrac{1}{n!}>2n\quad \textrm{\bigg(resp.  $RR_X(2)-\dfrac{c_X2^n}{(2n)!}=(n+1)^2-\dfrac{1}{n!}>2n$\bigg)}.
\end{align}
Same argument works for $\mathscr{M}_6$ as well since its Riemann-Roch expression is same as that of $K^3(T)$. That conludes the proof of the assertion.\QEDB\par
The upper bound for the degree has the following consequence on the secant lines of an embedding of a $K3$ surface. Even though it has nothing to do with our present purpose, it still might be of independent interest.

\begin{corollary}
\phantomsection\label{randomcor} Let $S$ be a $K3$ surface and $B$ be a very ample line bundle. Consider the projective embedding of $S$ in $\mathbb{P}^{h^0(B)-1}$ and the closed subvariety of $Gr(2,h^0(B))$ consisting of lines that intersect the $K3$ surface at a subscheme of length at least $2$. Then a general such line intersects $S$ at a subscheme of length $\leq 7$.
\end{corollary}
\noindent\textit{Proof.} Given the above conditions we construct a generically finite morphism $f$ from $X=S^{[2]}$ to $Gr(2,h^0(B))$. Given a point on $X$ we take the length $2$ subscheme it defines on $S$ and send it to the linear span of the length two subscheme inside $\mathbb{P}^{h^0(B)-1}$. Since a general such line does not lie on $S$, it intersects $S$ at finitely many points. So a general point in the image of the morphism $f$ has got finite fibers. Hence $f$ is a generically finite morphism.\\
\indent Also this morphism is given by the complete linear series $B-\delta$ where $2\delta$ is the class of the divisor in $S^{[2]}$ that parametrizes non-reduced subschemes of length $2$ on the $K3$ surface $S$. By Lemma \hyperref[ub]{3.4} and Remark \hyperref[hypcheck]{3.6}, we have $deg(f) \leq 23$ . If a line intersects $S$ at $k$ points then the line has $\binom{k}{2}$ preimages under the morphism $f$. Thus for a general line $\binom{k}{2}$   $\leq 23$ and hence $k \leq 7$. \QEDB\par

In the next Lemma, we will use the upper bound for $deg(\varphi_B)$ to put more restrictions on the image of $\varphi_L$ when it maps onto a variety of minimal degree.
\begin{lemma}\phantomsection\label{moreelimination}
Let $X$ be a hyperk\"ahler manifold of dimension $2n\geq 4$ for which $R_X\vert_{\mathbb{Z}\geq 0}$ is increasing, $R_X(\alpha_X)>2n$ and $RR_X(\alpha_X)\geq 4n$. Then for any ample and globally generated line bundle $B$ on $X$, $\varphi_B$ can never map $X$ onto a variety that is obtained by taking cones over a rational normal curve.
\end{lemma}
\noindent \textit{Proof.} Suppose $Y=Im(\varphi_B)$ is a variety of minimal degree that is obtained by taking cones over a rational normal curve and $d=deg(\varphi_B)$. Therefore $Y$ is singular in codimension two. Let $G$ be the inverse image of a general linear subspace of $R$ of codimension $1$ in $Y$. Notice that $G$ is irreducible by Bertini's theorem.\\
\indent We have $B^{2n-1}\cdot G = d$. Using Lemma \hyperref[rita]{1.10}, we deduce that $B$ can be written as the pullback of $deg(Y) \cdot R$. Thus, $G$ is ample and $d = deg(Y)^{2n-1}\cdot G^{2n} \geq deg(Y)^{2n-1}$. \\
\indent Now we use the fact that $deg(Y) = \left(RR_X(q_X(B)) - 2n \right)$ and $d<(2n)!$ (see Lemma \hyperref[ub]{3.4} and Remark \hyperref[hypcheck]{3.6}) that leads us to the following inequality,
\begin{align}
    \left(RR_X(q_X(B)) - 2n \right)^{2n-1} < (2n)!
\end{align} 
which is absurd. Indeed, by our assumption, $RR_X(q_X(B))-2n\geq RR_X(\alpha_X)-2n\geq 2n$. \QEDB\par 
Notice that the proof above also shows the following.
\begin{remark}\phantomsection\label{more..l}
Let $X$ be a hyperk\"ahler manifold of dimension $2n\geq 4$ for which $R_X\vert_{\mathbb{Z}\geq 0}$ is increasing and $R_X(\alpha_X)>2n$. Let $B$ be an ample and globally generated line bundle for which either $h^0(B)\geq 4n$, or $B^{2n}<\alpha_X^n(h^0(B)-2n)^{2n}$. Then $\varphi_B$ can never map $X$ onto a variety that is obtained by taking cones over a smooth rational normal curve.
\end{remark}
Combining Lemmas \hyperref[elimination]{3.1}, \hyperref[moreelimination]{3.8} and Remarks \hyperref[hypcheck]{3.6} and \hyperref[more..l]{3.9} we get the following.
\begin{proposition}\phantomsection\label{hkvmd}
Let $X$ be a hyperk\"ahler manifold of dimension $2n\geq 4$ and let $B$ be an ample and globally generated line bundle on $X$. If $X$ is deformation equivalent to either $K^n(T)$ or $\mathscr{M}_6$ then the morphism $\varphi_B$ given by the complete linear series $|B|$ will never map $X$ onto a variety of minimal degree. If $X$ is of type $K3^{[n]}$ and $\varphi_B$ maps $X$ onto a variety $Y$ of minimal degree then either, 
\begin{itemize}
    \item[(1)] $X$ is of type $K3^{[2]}$, $q_X(B)=2$, $deg(\varphi_B)=6$ and $Y$ is a quadric hypersurface (possibly singular) inside $\mathbb{P}^5$ which can not be obtained by taking cones over any rational normal scroll, or
    \item[(2)] $X$ is of type $K3^{[3]}$, $q_X(B)=2$, $deg(\varphi_B)=30$ and $Y$ is a variety embedded inside $\mathbb{P}^9$ that is obtained by taking cones over the Veronese embedding of $\mathbb{P}^2$ inside $\mathbb{P}^5$.
\end{itemize}
\end{proposition} 
\noindent\textit{Proof.} To start with, note that $\alpha_X\geq 2$. Suppose $X$ is of type $K^n(T)$ or $\mathscr{M}_6$. By Riemann-Roch, we get $h^0(B)=RR_X(q_X(B))\geq RR_X(2)>2n+4$. Consequently, by Lemma \hyperref[elimination]{3.1}, $Im(\varphi_B)$ can only be a variety that is obtained by taking cones over a rational normal curve. But that is impossible by Lemma \hyperref[moreelimination]{3.8} since $RR_X(2)>4n$.\\
\indent Now, assume $X$ is of type $K3^{[n]}$. We can argue exactly like the previous paragraph to conclude that $\varphi_B$ will never map $X$ onto a variety of minimal degree if $n\geq 5$.\\
\indent  We deal with the case $n=2$, $3$ and $4$ separately and we will use Lemma \hyperref[elimination]{3.1}. Note that $q_X(B)$ is even, say $q_X(B)=2k$ for some positive integer $k$.\\
\indent Suppose $n=2$. Note that $h^0(B)>8$ if $q(B)\geq 4$ and $B^{2n}<2^2(h^0(B)-4)^4$ if $q(B)=2$. Consequently by Remark \hyperref[more..l]{3.9}, $Y$ can not be obtained by taking cones over rational normal curve. The equation $RR_X(2k)=2n+2$ has only one positive even integer solution $k=1$ in which case $q_X(B)=2$, $deg(\varphi_B)=6$ and $Y$ is a quadric hypersurface in $\mathbb{P}^5$. $RR_X(2k)=2n+4$ has no integer solution.\\
\indent Suppose $n=3$. Argument similar to that of the case $n=2$ yields that $Y$ can not be obtained by taking cones over rational normal curve. Moreover, $RR_X(2k)=2n+2$ has no solution and $RR_X(2k)=2n+4$ has only one positive integer solution $k=1$ in which case $q_X(B) = 2$, $deg(\varphi_B) = 30$ and $Y$ is a variety of minimal degree in $\mathbb{P}^9$ obtained by taking veronese embedding of $\mathbb{P}^2$ inside $\mathbb{P}^5$.\\
\indent Similar argument shows that $\varphi_B$ can not map $K3^{[4]}$ onto a variety of minimal degree.\QEDB\par Now we are ready to give the proof of Theorem \hyperref[ob]{B}. 
\begin{theorem}\phantomsection\label{pthmb}
Let $X$ be a projective hyperk\"ahler variety of dimension $2n\geq 4$ that is deformation equivalent to $K3^{[n]}$, $K^n(T)$ or $\mathscr{M}_6$. Let $B$ be an ample and globally generated line bundle on $X$. Then the following happens;
\begin{itemize}
    \item[(1)] $B^{\otimes l}$ is projectively normal for $l\geq 2n$.
    \item[(2)] $B^{\otimes 2n-1}$ is projectively normal unless:
    \begin{itemize}
        \item [(a)] $n=2$, $X=K3^{[2]}$ and $\varphi_B$ maps $X$ onto a quadric (possibly singular) inside $\mathbb{P}^5$. In this case $q_X(B)=2$, $deg(\varphi_B)=6$, or
        \item[(b)] $n=3$, $X=K3^{[3]}$ and $\varphi_B$ maps $X$ onto a variety of minimal degree inside $\mathbb{P}^9$ which is obtained by taking cones over the Veronese embedding of $\mathbb{P}^5$ inside $\mathbb{P}^5$. In this case $q_X(B)=2$, $deg(\varphi_B)=30$.
    \end{itemize}
\end{itemize}
Hence if $X$ is as above and $B$ does not satisfy cases $2(a)$ or $2(b)$ then a general curve section of $|B|$ is non-hyperelliptic.
\end{theorem}
\noindent\textit{Proof.} To prove (1) we simply notice that $h^0(B) \geq 2n+2$ by the Riemann-Roch formula on $X$. The assertion follows by Lemma \hyperref[old]{2.1}. (2) follows directly by Corollary \hyperref[cor]{1.8}, and Proposition \hyperref[hkvmd]{3.10}. The statement on non-hyperellipticity of a general curve section $C$ follows from the fact that $B^{d-1}|_C = K_C$ by adjunction and that a very ample line bundle restricts to a very ample line bundle on a closed immersion.  \QEDB\par 
We finish by the following example of Debarre (see \cite{De}). It shows the existence of an ample and globally generated line bundle on a hyperk\"ahler variety of $K3^{[2]}$ type that induces a 6-1 map onto a variety of minimal degree by its complete linear series.\par 
\noindent {\bf Example 3.12}\phantomsection\label{ex} Let $(S,L)$ be a polarized $K3$ surface with $Pic(S)=\mathbb{Z}L$ and $L^2=4$. Then $L$ is very ample and consequently we get a morphism $\phi: S^{[2]}\rightarrow Gr(2,4)$ to the Grassmannian.\\
\indent Now, $L$ induces a line bundle $L_2$ on $S^{[2]}$ and it is known that $Pic(S^{[2]})=\mathbb{Z}L_2\oplus\mathbb{Z}\delta$. Moreover, the pullback of the Pl\"ucker line bundle on the Grassmannian has class $L_2-\delta$ on $S^{[2]}$. Therefore, if $(S,L)$ is general then it contains no line and consequently $\phi$ will be finite of degree $\binom{4}{2} =6$.\QEDB

\bibliographystyle{plain}

\end{document}